\documentclass[twocolumn,floats,groupedaddress,aps]{revtex4}
\usepackage{times,amsthm,amsmath,amsfonts,latexsym}
\usepackage{pstricks,pst-node,graphics,psfrag}
\makeatletter
\renewcommand{\thesection}{\arabic{section}}
\renewcommand{\p@subsection}{}
\renewcommand{\thesubsection}{\arabic{section}.\arabic{subsection}}
\renewcommand{\p@subsubsection}{}

\makeatother

\newtheorem{theorem}{Theorem}

\newtheorem{problem}[theorem]{Problem}
\newtheorem{corollary}[theorem]{Corollary}
\newtheorem{conjecture}[theorem]{Conjecture}

\theoremstyle{remark}
\newtheorem*{remark}{Remark}

\DeclareMathOperator{\Tr}{Tr}
\DeclareMathOperator{\Vol}{Vol}

\newcommand{\GL}{\mathrm{GL}} \newcommand{\GUE}{\mathrm{GUE}}
\newcommand{\LP}{\mathrm{LP}} \newcommand{\QM}{\mathrm{QM}}
\newcommand{\RSK}{\mathrm{RSK}} \newcommand{\SO}{\mathrm{SO}}
\newcommand{\SU}{\mathrm{SU}} \newcommand{\U}{\mathrm{U}}

\newcommand{\gl}{\mathfrak{gl}}
\newcommand{\mg}{\mathfrak{g}}
\newcommand{\su}{\mathfrak{su}}
\newcommand{\so}{\mathfrak{so}}
\renewcommand{\sp}{\mathfrak{sp}}
\renewcommand{\sl}{\mathfrak{sl}}
\renewcommand{\u}{\mathfrak{u}}
\newcommand{\mC}{\mathfrak{C}}
\newcommand{\mh}{\mathfrak{h}}

\renewcommand{\hat}{\widehat} \renewcommand{\tilde}{\widetilde}

\newcommand{\tA}{\tilde{A}} \newcommand{\tB}{\tilde{B}}
\newcommand{\tJ}{\tilde{J}} \newcommand{\tM}{\tilde{M}}

\newcommand{\hdelta}{{\hat{\delta}}} \newcommand{\hlambda}{{\hat{\lambda}}}
\newcommand{\hmu}{{\hat{\mu}}} \newcommand{\hp}{\hat{p}}

\newcommand{\cN}{\mathcal{N}} \newcommand{\cH}{\mathcal{H}}
\newcommand{\cM}{\mathcal{M}}

\newcommand{\C}{\mathbb{C}} \newcommand{\Z}{\mathbb{Z}}
\newcommand{\R}{\mathbb{R}}

\newcommand{\ket}[1]{|#1\rangle}

\newcommand{\bracket}[1]{\langle#1\rangle}
\newcommand{\up}{\uparrow}
\newcommand{\down}{\downarrow}
\newcommand{\aket}[1]{|\!#1\,\rangle}
\newcommand{\abra}[1]{\langle\,#1\!|}
\newcommand{\st}{\:|\:}

\newcommand{\eatline}{\vspace{-\baselineskip}}
\newcommand{\ie}{\textit{i.e.}}

\newcommand{\eq}[2]{\begin{equation}\label{#1}#2\end{equation}}
\renewcommand{\matrix}[2]{\left(\begin{array}{#1}#2\end{array}\right)}
\renewcommand{\tensor}{\otimes}

\renewcommand{\figurename}{Figure}
\newenvironment{fullfigure}[2]
    {\begin{figure}[htb]\begin{center}\def\ffa{#1}\def\ffb{#2}}
    {\vspace{\baselineskip}\caption{\ffb.}\label{\ffa}\end{center}\end{figure}}
\psset{linewidth=.5pt,dash=4pt 4pt,unit=.25in}
\newgray{gray40}{.6}
\SpecialCoor

\begin{document}
\title{Random words, quantum statistics, central limits, random matrices}
\author{Greg Kuperberg}
\thanks{Supported by NSF grant DMS \#9704125 and by a
    Sloan Foundation Research Fellowship}
\affiliation{UC Davis; {\tt greg@math.ucdavis.edu}}
\begin{abstract}
Recently Tracy and Widom conjectured \cite{TW:monotone} and Johansson proved
\cite{Johansson:plancherel} that the expected shape $\lambda$ of the
semi-standard tableau produced by a random word in $k$ letters is
asymptotically the spectrum of a random traceless $k \times k$ GUE matrix.  In
this article we give two arguments for this fact. In the first argument, we
realize the random matrix itself as a quantum random variable on the space of
random words, if this space is viewed as a quantum state space. In the second
argument, we show that the distribution of  $\lambda$ is asymptotically given
by the usual local limit theorem, but the resulting Gaussian is disguised by an
extra polynomial weight and by reflecting walls.  Both arguments more generally
apply to an arbitrary finite-dimensional representation $V$ of an arbitrary
simple Lie algebra $\mg$.  In the original question, $V$ is the defining
representation of $\mg = \su(k)$.
\end{abstract}
\maketitle

What is the longest weakly increasing subsequence of a long, random string of
letters?  In the previous sentence, one such longest subsequence is
``AEEEEEEEFLNNOSTTT''.  In randomly chosen English text, the longest
subsequences are dominated by the letter 'E', since this letter is the most
common one.  This implies that the length of the longest subsequence has a
Gaussian distribution.   But if the letters in the string are independent with
the uniform distribution, a longest subsequence will use all of them roughly
equally.  In this case Tracy and Widom established a non-Gaussian distribution
for the length of a longest subsequence \cite{TW:monotone,ITW:toeplitz}. Their
result was motivated by recent progress in the study of the longest increasing
subsequence of a random permutation, in particular the relations among longest
subsequences, random matrices, and representation theory
\cite{AD:patience,BDJ:longest,BDJ:second,TW:painleve,BO:z-measures,BOO:plancherel,%
BR:algebraic,Okounkov:permutations}.

Tracy and Widom conjectured a generalization which
was proved by Johansson \cite[Th. 1.6]{Johansson:plancherel}:

\begin{theorem}[Johansson] The distribution of the shape of a random word as
given by the Robinson-Schensted-Knuth (RSK) algorithm converges locally to
the distribution of the spectrum of a random traceless $k \times k$ GUE matrix.
\label{th:main}
\end{theorem}

It is a generalization because the first row of the RSK shape is the length of
the longest weakly increasing subsequence. ``Traceless GUE'' refers to the
traceless Gaussian unitary ensemble, defined up to normalization as the
Gaussian measure on traceless $k \times k$ Hermitian matrices which is
invariant under conjugation by unitary matrices.

In this article we give two arguments for Theorem~\ref{th:main}. The first
argument (Section~\ref{s:quantum}) is based on quantum statistics:  it identifies the
random matrix itself as a quantum random variable on the space of random words
viewed as a quantum state space.  The GUE ensemble then appears in the limit by
a quantum central limit theorem.  The second argument (Section~\ref{s:local})
is based on classical statistics:  it identifies the density formula
\eq{e:mehta}{C\prod_{a \le b}
    (\lambda_a - \lambda_b)^2 e^{-\sum_a \lambda_a^2} d\lambda}
for the distribution of the spectrum $\lambda$ of a GUE matrix
\cite{Mehta:matrices} as a disguised classical central limit.  (Here $C$ is a
constant that depends on $k$ but not $\lambda$.) The classical argument is
rigorous and it establishes a precise estimate. The quantum argument can be
read rigorously or non-rigorously, depending on whether the reader accepts
Conjecture~\ref{c:qcentral}; either way it is less precise than the statement
of Theorem~\ref{th:main}.  Non-rigorously, it establishes convergence in
distribution.  Rigorously it establishes convergence of certain moments, but
not enough moments to imply convergence in distribution. Nonetheless we prefer
the quantum argument since it is less traditional.  (But see Biane
\cite{Biane:bernoulli,Biane:dual,Biane:tensorielles} for closely
related results.)

In both arguments, it is important to identify the vector space of traceless
Hermitian matrices with the Lie algebra $\su(k)$ and an alphabet with $k$
letters with the standard basis of the defining representation $V = \C^k$. 
Both arguments then generalize to an arbitrary finite-dimensional unitary
representation of $V$ of a compact simple Lie algebra $\mg$.  The conclusion is
a relation between random words in a weight basis of $V$ and a natural Gaussian
measure on $\mg^*$, the vector space dual of $\mg$.

\acknowledgments

The author would like to thank Greg Lawler, Marc van Leeuwen, Bruno
Nachtergaele, Philippe Biane, and Marc Rieffel for informative discussions, and
especially Craig Tracy for introducing him to the topics discussed in this
article.

\section{Quantum statistics}
\label{s:quantum}

In this section we will express certain classical random variables in terms of
simpler quantum random variables. The main object in our argument was also
considered from the converse view by Biane \cite{Biane:dual,Biane:bernoulli}.

We refer the reader to Sakurai \cite[\S3]{Sakurai:modern} for basic notions of
quantum statistics, in particular \emph{mixed states}, which are also commonly
called density matrices or density operators.  In the context of operator
algebras, mixed states are called \emph{states} \cite{KR:vol1} or \emph{normal
states} \cite{KR:vol2}, depending on the desired strength of the formalism.

The \emph{RSK algorithm} is (in one version) a function that takes as input a word of
length $N$ in the alphabet $[k] = \{1,...,k\}$ and produces as output a pair of
tableaux $(P,Q)$ of shape $\lambda$, where
$$\lambda = (\lambda_1, \lambda_2,\ldots, \lambda_k) \vdash N$$
is a partition of $N$ into non-increasing, non-negative integers
\cite[\S7.11]{Stanley:vol2}.  The partition is considered synonymous with its \emph{Young
diagram}, meaning its horizontal histogram. The tableau $P$ is semi-standard
and is called the  \emph{insertion tableau}, while the tableau $Q$ is standard
and is called the \emph{recording tableau}. Given the uniform distribution on the
set of words $[k]^N$, we can view the shape $\lambda$ as a random variable
$\lambda_\RSK$. Finally, given a partition $\lambda$, it will sometimes be
convenient to subtract the mean from each part to form  a ``partition of 0'':
$$\hlambda = (\lambda_1-\frac{N}k,\lambda_2-\frac{N}k,\ldots,
    \lambda_k-\frac{N}k).$$

We do not need the precise definition of the RSK algorithm in this section,
merely one of its important properties: It is a combinatorial model for the
direct sum decomposition of the representation $V^{\tensor N}$ of the
Lie algebra $\u(k)$ (or the Lie group $\U(k)$ or $\GL(k,\C)$),
where $V = \C^k$ is the defining representation \cite[\S A2]{Stanley:vol2}. This
representation decomposes as
\eq{e:decomp}{V^{\tensor N} \cong \bigoplus_{\lambda\vdash N} R_\lambda
    \tensor V_\lambda,}
where $V_\hlambda$ is the irreducible representation of $\u(k)$ of shape
$\lambda$ and $R_\lambda$ is the irreducible representation of the
symmetric group $S_N$ of shape $\lambda$. For any given $\lambda =
\lambda_\RSK$, the set of associated insertion tableaux $P$ indexes a basis of
$V_\lambda$, while the set of recording tableaux $Q$ indexes a basis of
$R_\lambda$.  In particular,
$$\dim R_\lambda \tensor V_\lambda = n_\lambda,$$
where $n_\lambda$ is the number of words that have shape $\lambda =
\lambda_\RSK$.  Finally, as a representation of the Lie subalgebra
$\su(k)$, $V_\lambda$ is unchanged if we add a constant to each
component of $\lambda$.  The convention is to call it $V_\hlambda$,
the representation of highest weight is $\hlambda$.

We can view the vector space $V^{\tensor N}$ as a quantum state space $\cH$ of
some quantum system $Q$ with $[k]^N$ as an orthonormal basis.  The
maximum-entropy state (or \emph{tracial} state) $\rho$ of $Q$ is then realized
by the uniform distribution on $[k]^N$, as well as by the uniform distribution
on any other orthonormal basis.  At the same time, an arbitrary orthogonal
direct sum decomposition
$$\cH \cong \cH_1 \oplus \cH_2 \oplus \ldots \oplus \cH_t$$
of $\cH$ can be interpreted as a random variable taking values in the set of
summands.  Relative to the state $\rho$, the probability of a given summand
$\cH_i$ is the ratio $(\dim \cH_i)/(\dim \cH)$.  In particular,  the direct sum
decomposition in equation~\eqref{e:decomp} expresses a random variable
$\lambda_\QM = \hlambda$.  The previous paragraph tells us that
$\lambda_\QM \doteq \hlambda_\RSK$, meaning that they have the same 
distribution.

\subsection{The case $k=2$ and spin 1/2 particles}
\label{s:spin}

As a concrete example, we consider the physically realizable case $k=2$.  In
this case $V$ is the familiar state space of a spin $1/2$ particle, and the
action of $\SU(2)$ is the projective action of the spatial rotation group
$\SO(3)$.   We will use the alphabet $\{\up,\down\}$ rather than $\{1,2\}$ as a
basis of $V$.  The space $V$ admits angular momentum operators $J_x$, $J_y$,
and $J_z$ which satisfy the commutation relations
\eq{e:su2}{[J_x,J_y] = iJ_z \qquad [J_y,J_z] = iJ_x \qquad [J_z,J_x] = iJ_y.}
The operators $J_x$, $J_y$, and $J_z$ are a basis of $i \cdot \su(2)$, by which
we mean the image of $\su(2)$ in $\sl(2,\C)$ under multiplication by
$i$. Thus these are just the usual commutation relations in the Lie
algebra $\su(2)$ up to a factor of $i$.  The tracial state on $V$ is the
mixture
\eq{e:mix}{\rho = \frac{\aket\up\abra\up + \aket\down\abra\down}{2}.} 
Note that probabilistic mixtures of states should not be confused with quantum
superpositions. A superposition of $\aket\up$ and $\aket\down$ is another
vector in $V$ and cannot be invariant under rotations. By contrast the mixed
state $\rho$ is $\SU(2)$-invariant.

The vector space $V^{\tensor N}$ is then the state space of $N$ such particles.
The $t$th particle has angular momentum operators $J_x^{(t)}$, $J_y^{(t)}$, and
$J_z^{(t)}$. By equation \eqref{e:mix} and the isotropy of $\rho$, each of
these operators is a centered Bernoulli random variable with equally likely
values $\frac12$ and $-\frac12$.  Since the three operators for any fixed $t$
do not commute, the corresponding random variables cannot be simultaneously
observed.  The sums of these operators form the total angular momentum,
\eq{e:sum}{J_\alpha =
    J_\alpha^{(1)} + J_\alpha^{(2)} + \ldots + J_\alpha^{(N)},}
for each $\alpha \in \{x,y,z\}$.  Each operator $J_\alpha$ is a centered
binomial random variable because the terms are independent commuting Bernoulli
variables.  The three operators $J_x$, $J_y$, and $J_z$ do not commute either,
but rather satisfy the same commutation relations in equation~\eqref{e:su2},
since they express the  natural (diagonal) action of $\su(2)$ on $V^{\tensor
N}$.   Finally, the total angular momentum
$$J^2 = J_x^2 + J_y^2 + J_z^2$$
is diagonalized by the direct sum decomposition in equation~\eqref{e:decomp}.
In a summand with weight $\hlambda$, its eigenvalue is
$$J^2\ket\psi = \hlambda_1(\hlambda_1 + 1)\ket\psi.$$
Thus if $L = (L_1,L_2)$ is the shape-valued operator that measures
$\lambda_\QM$, it is related to $J^2$ by
$$J^2 = L_1(L_1+1).$$

If we define the scaled angular momentum operators
$$\tJ_\alpha = \frac{J_\alpha}{\sqrt{N}}$$
with $\alpha \in \{x,y,z\}$, then by the above reasoning, these three operators
become commuting Gaussian random variables in the limit $N \to \infty$.  By
rotational symmetry they are also independent and identically distributed
(i.i.d.). 
This behavior of the total angular momentum of independent random spins can be
witnessed in nuclear magnetic resonance experiments, among other places.  Using
these operators we may form a matrix
$$\tM = \matrix{cc}{\tJ_z & \tJ_x + i\tJ_y \\ \tJ_x - i\tJ_y & -\tJ_z}.$$
In the limit $N \to \infty$, $\tM$ becomes a traceless GUE matrix!  (The
normalization is also consistent with Mehta \cite{Mehta:matrices}.)  For finite
$N$, the determinant of $\tM$ must be interpreted carefully because its entries
do not commute.  If we define it by averaging over orderings of the entries,
$$\det \tM = \frac12 (\tM_{11}\tM_{22} + \tM_{22}\tM_{11}
    - \tM_{12}\tM_{21} - \tM_{21} \tM_{12}),$$
then it turns out that
$$\det \tM = -\tJ^2.$$
It follows that
$$\lim_{N \to \infty} \frac{\lambda_\QM}{\sqrt{N}} \doteq \lambda_\GUE,$$
where $\lambda_\GUE$ is a random variable representing the spectrum $\lambda$
of a traceless GUE matrix. This is precisely Theorem~\ref{th:main} for $k=2$.

\subsection{The general case}

The argument in Section~\ref{s:spin} generalizes with only natural changes to
all values of $k$.  The defining representation $V$ of $\su(k)$ has a basis of
states
$$\ket1,\ket2,\ldots,\ket{k}.$$
The elements of $i\cdot\su(k)$ may be viewed as generalized angular momentum
operators.  We define two matrices of operators $A$ and $B$ whose entries
linearly span $i \cdot \su(k)$.  (Note that the diagonal entries are not
linearly independent.)  Let $E_{ab} \in \sl(k,\C)$ be the elementary matrix
whose non-zero entry is
$$(E_{ab})_{ab} = 1.$$
Then the entries of $A$ and $B$ are
\begin{align*}
A_{ab} &= \frac12 (E_{ab} + E_{ba}) &
B_{ab} &= \frac{i}2 (E_{ba} - E_{ab}) \qquad (a \ne b) \\
A_{aa} &= E_{aa} - \frac1{k}I & B_{aa} &= 0.
\end{align*}
Let $M$ be the matrix of operators
$$M_{ab} = A_{ab} + iB_{ab}.$$

Each entry $A_{ab}$ and $B_{ab}$ is a real-valued measurement operator.
Relative to the tracial state on $V$, $A_{ab}$ for $a \ne b$ takes each of
the values $1/2$ and $-1/2$ with probability $1/k$ and otherwise has the value
$0$.  The same is true of $B_{ab}$.  The measurement $A_{aa}$ takes the value
$(k-1)/k$ with probability $1/k$ and otherwise takes the value $-1/k$.  The
operator $M_{ab}$ may appear to be a complex-valued measurement whose real and
imaginary parts have these distributions, but this is not quite true. An
operator can only be interpreted as a complex-valued measurement if it is a
normal operator, defined as an operator that commutes with its adjoint, or
equivalently an operator whose self-adjoint and anti-self-adjoint parts
commute.  When $a \ne b$, $M_{ab}$ is not a normal operator; it represents a
complex random variable whose real and imaginary parts are not simultaneously
observable.

For words in $[k]^N$, we consider the standard (additive) action of $\su(k)$ on
$V^{\tensor N}$.  Equation~\eqref{e:decomp} gives us a shape-valued operator
$L$ which has the eigenvalue $\lambda_{\QM} = \hlambda$ on the summand
$V_\lambda \tensor R_\lambda$. The operator $L$ can be realized algebraically
using the characteristic polynomial of $M$, thought of as a polynomial-valued
operator:
\begin{align} C(x) &= \det\:(xI - M) \nonumber \\
&= (x - L_1)(x - L_2)\ldots(x - L_k) + c(x,L), \label{e:det}
\end{align}
where $c(x,L)$ is a polynomial of total degree at most $k-1$. As before, each
term of the determinant is defined by averaging over the $k!$ orderings of its
factors.  The left side of equation~\eqref{e:det} is a disguised version of the
composition $\psi \circ \pi$ in the proof of Harish-Chandra's Theorem given by
Humphreys \cite[\S23.3]{Humphreys:gtm}, while the leading term on the right side is
a disguised version of the map $\theta$.  Here we are applying both maps to the
coefficients of the ordinary characteristic polynomial of an element in
$\su(k)$; in the context of Humphreys, each coefficient is a particular
$\SU(k)$-invariant polynomial on $\su(k)$.  As Humphreys explains, the maps
$\psi \circ \pi$ and $\theta$ agree in the top degree, which is exactly what
equation~\eqref{e:det} asserts.  (See Okounkov and Olshanski
\cite{OO:shifted} for an analysis of the correction term $c(x,L)$.)

Each coefficient of $C(x)$ lies in the center of $U(\su(k))$ and is a natural
generalization of the Casimir operator $J^2$ in the case $k=2$.  The
coefficients are sometimes called elementary generalized Casimir operators.

Assuming the tracial state on $V^{\tensor N}$, each measurement $A_{ab}$ and
$B_{ab}$ is a sum of bounded, centered i.i.d. random variables.  If we define
$$\tM = \tA + i\tB = \sqrt{\frac{k}{2N}} M,$$
then the entries commute in the limit $N \to \infty$ and $\tM$ becomes a
traceless GUE matrix with standard normalization.  The term
$c(x,L)$ in equation~\eqref{e:det} also disappears in this limit
because its degree is too low.  The equation thus tells us that
$$\lim_{N \to \infty} \sqrt{\frac{k}{2N}}\lambda_\QM \doteq \lambda_\GUE.$$

\subsection{What did we prove?}

One important step in the argument of this section is not completely rigorous.
Unquestionably each scaled angular momentum operator $\tJ_\alpha$, $\tA_{ab}$,
or $\tB_{ab}$ converges to a Gaussian random variable by the classical central
limit theorem. But we do not know that a polynomial in these variables, for
example $\tJ^2$ or $\tJ_x \tJ_y \tJ_z$, converges in distribution to the
corresponding polynomial of Gaussian variables. We cannot appeal to the
classical multivariate central limit theorem for non-commuting variables, even
if they do commute in the limit.  There are also several quantum limit theorems
in the literature; one of the most general ones is due to Goderis, Verbeure,
and Vets \cite{GVV:central}. But these results are apparently not sufficiently
strong either.

As a stop-gap we will conjecture the quantum central limit theorem that we
need, and we will prove a weak version of the conjecture.  The conjecture is
naturally stated in terms of $C^*$-algebras and von Neumann algebras
\cite{KR:vol1,KR:vol2}, which provide a rigorous language for quantum
statistics.  In this language, a non-commutative probability space is defined
as a von Neumann algebra $\cM$ with a normal state $\rho$. A \emph{state} is
defined as a dual vector on $\cM$, continuous in the norm topology, with the
interpretation that for a self-adjoint element $A$, $\rho(A)$ is the expected
value of the random variable given by $A$; the state is \emph{normal} if it is
continuous in the weak topology as well.  The reader who is
uninterested in operator algebras can take
$$\cM = \cM_k = \gl(k,\C),$$
the vector space of $k \times k$ matrices.  A normal state $\rho$ on $\cM_k$ is
any dual vector whose matrix is Hermitian and positive semi-definite and has
trace 1; this is exactly the definition in physics of a finite density matrix
\cite[\S3]{Sakurai:modern}. In particular, the tracial state is defined by
$$\rho(A) = \frac1k \Tr A.$$

Given two quantum systems with von Neumann algebras $\cM$ and $\cN$, the joint
system has the algebra $\cM \tensor \cN$, using the tensor product in the
category of von Neumann algebras.  Two normal states $\rho$ and $\omega$ on
$\cM$ and $\cN$ form an independent joint state $\rho \tensor \omega$.  Given
self-adjoint operators $A \in \cM$ and $B \in \cN$, the operators $A \tensor I$
and $I \tensor B$ represent independent measurements in the joint system. If $A
\in \cM$ is self-adjoint, let
$$A^{(t)} = I^{\tensor t-1} \tensor A \tensor I^{\tensor N-t-1}
\in \cM^{\tensor N},$$
and let
$$\tA = \frac1{\sqrt{N}} \sum_{t=1}^N A^{(t)}.$$
Thus $\tA$ expresses the scaled sum of $N$ i.i.d. random variables.

In formulating a multivariate quantum central limit theorem, three issues arise
because of non-commutativity.  First, the theorem must be a statement about the
distribution of non-commutative polynomials $p \in \C\bracket{A_1,\ldots,A_k}$,
but a Gaussian central limit would describe the distribution of commuting
variables. Second, a general polynomial expression $p(A_1,\ldots,A_k)$ need
not be a self-adjoint operator, even if the variables are self-adjoint.  Thus
we will assume that $p$ is a \emph{self-adjoint polynomial}, meaning that it is
invariant under the anti-linear anti-involution
$$*:\C\bracket{A_1,\ldots,A_k} \to \C\bracket{A_1,\ldots,A_k}$$
that conjugates each coefficient and reverses the order of each term.  Third,
if we define the covariance matrix  of the variables $A_1,\ldots,A_k$ as
$$\kappa_{ab} = \rho(A_a A_b),$$
it may not be symmetric.  When this happens the behavior of the central limit
is genuinely different from the classical case; it has been studied in
Reference~\citealp{GVV:central}. But in the case that we need (namely, when
$\rho$ is tracial), the covariance matrix is symmetric. In this case we expect
that the limiting distribution of $p$ only depends on its commutative image
$\hp \in \C[A_1,\ldots,A_k]$, thereby resolving the first issue.

\begin{conjecture} Let $(\cM,\rho)$ be a quantum probability space, and let
$A_1,\ldots,A_k$ be self-adjoint elements with mean 0 and a symmetric
covariance matrix. Let $p \in \C\bracket{A_1,\ldots,A_k}$ be a self-adjoint
non-commutative polynomial in $k$ variables.  Then
$$\lim_{N \to \infty} p(\tA_1,\ldots,\tA_k) \doteq \hp(X_1,\ldots,X_k),$$
where $X_1,\ldots,X_k$ are classical Gaussian random variables with
covariance matrix
$$E[X_a X_b] = \rho(A_a A_b).$$
\label{c:qcentral}
\eatline\end{conjecture}

If we let $p$ be a coefficient of the polynomial $C(x)$ from
Section~\ref{s:quantum}, Conjecture~\ref{c:qcentral} then implies
Theorem~\ref{th:main}.  It may be possible to be reverse this reasoning and use
Theorem~\ref{th:main} to prove Conjecture~\ref{c:qcentral} for arbitrary $p$,
at least when $\cM$ is a matrix algebra and $\rho$ is the tracial state.   But
a satisfactory proof would hold for arbitrary quantum probability spaces.

Our Theorem~\ref{th:qcentral} below establishes convergence of moments in the
context of Conjecture~\ref{c:qcentral}.   Since this theorem is almost entirely
algebraic, we do not need the full structure of a von Neumann algebra. Rather
we let $\cM$ be an arbitrary $*$-algebra, meaning a unital ring over the
complex numbers with an anti-linear anti-involution $*$.  A state $\rho$ on a
$*$-algebra is a $*$-invariant dual vector such that $\rho(I) = 1$ and
$\rho(A^2) \ge 0$ for every self-adjoint $A \in \cM$.  Finally
the $n$th moment
$$\gamma_n(A) = \rho(A^n)$$
is defined whether or not $A$ is self-adjoint. (Recall that that a
non-self-adjoint operator may be written as $A+iB$, where $A$ and $B$ are
self-adjoint.  Consequently it may be interepreted as a complex-valued quantum
random variable whose real and imaginary parts are not simultaneously
observable.  Our definition of moments is consistent with this
interpretation.)

\begin{theorem} Suppose that $\cM$ a $*$-algebra with a state $\rho$. Let
$A_1,\ldots,A_k$ be self-adjoint elements with mean 0. Suppose that 
for all $a$ and $b$,
$$\rho(A_a A_b) = \rho(A_b A_a).$$
If $p \in \C\bracket{A_1,\ldots,A_k}$ is a non-commutative
polynomial, then
$$\lim_{N \to \infty} \gamma_n(p(\tA_1,\ldots,\tA_k)) =
\gamma_n(\hp(X_1,\ldots,X_k)),$$
where $X_1,\ldots,X_k$ are classical centered Gaussian random variables with
covariance matrix 
$$E[X_a X_b] = \rho(A_a A_b),$$
and $\gamma_n(A)$ is the $n$th moment of $A$.
\label{th:qcentral}
\end{theorem}
\begin{proof} Since the assertion is claimed for every polynomial, it suffices
to prove that the expectation
$$\rho^{\tensor N}(p(\tA_1,\ldots,\tA_k))$$
converges.  To show convergence of expectation we may let $p$ be a monomial. 
Indeed the monomial
$$p(A_1,\ldots,A_k) = A_1A_2\ldots A_k$$
will do, since some of the factors may be equal.

Expanding the expression
$$\gamma = \rho^{\tensor N}(\tA_1,\ldots,\tA_k)$$
using the definition of $\tA_a$, it has a term for
each function $\phi$ from $[k]$ to $[N]$:
$$\gamma = N^{-k/2} \sum_{\phi:[k] \to [N]} \rho^{\tensor N}(\prod_a A_a^{(\phi(a))}).$$
Let
$$\gamma_\phi = N^{-k/2} \rho^{\tensor N}(\prod_a A_a^{(\phi(a))})$$
be an individual term in this expansion. Since we are computing the expectation
with respect to the  product state, we can arrange the factors with respect to
$[N]$ rather than $[k]$:
$$\gamma_\phi = N^{-k/2} \prod_{t\in [N]} \rho(\prod_{\phi(a) = t} A_a^{(t)}).$$
In this form it is clear that
$$\gamma_\phi = 0$$
if there is a $t$ such that $\phi^{-1}(t)$ has one element. At the same
time, if $S$ is the set of those functions $\phi$ whose images have fewer than
$k/2$ elements, then
$$\lim_{N \to \infty} \sum_{\phi \in S} \gamma_\phi = 0$$
because
$$|S| = o(N^{k/2}).$$
In other words, there are sufficiently few such functions $\phi$ that they are
negligible in the limit. What is left is the set of functions that are exactly
2-to-1, which only exist when $k$ is even.  Thus if $M$ is the set of perfect
matchings of $[k]$, then
\begin{align*}
\gamma =&\  N^{-k/2} \binom{N}{k/2}(k-1)!!(k/2)! \sum_{m \in M}
    \prod_{(a,b) \in m} \rho(A_a A_b) \\
& + o(1) \\
=&\  (k-1)!!\sum_{m \in M}  \prod_{(a,b) \in m} \rho(A_a A_b) + o(1),
\end{align*}
where
$$(k-1)!! = (k-1)(k-3)\cdots 5 \cdot 3 \cdot 1 $$
is the odd factorial function.  In the limit $\gamma$ exactly matches the
corresponding expectation
$$E[X_1\ldots X_k] = (k-1)!!\sum_{m \in M} \prod_{(a,b) \in m} E[X_a X_b]$$
of the classical variables $X_1,\ldots,X_k$.
\end{proof}

In relation to Theorem~\ref{th:main}, Theorem~\ref{th:qcentral}
says that if
$$C_\lambda(x) = (x-\lambda_1)(x-\lambda_2)\ldots(x-\lambda_k)$$
for a partition $\lambda$, then the moments of $C_{\hlambda_\RSK}(x)$ converge
to the moments of $C_{\lambda_\GUE}(x)$.  In other words, the ``moments of the
moments'' of $\lambda_\RSK$ converge after scaling to those of $\lambda_\GUE$. 
Unfortunately, when $k > 2$ the tail of the distribition of
$C_{\lambda_\GUE}(x)$ is too thick for convergence of moments to imply
convergence in distribution.

\section{Local limits}
\label{s:local}

The argument of this section and its generalization below
(Theorem~\ref{th:general}) are very similar to a result of Biane
\cite{Biane:tensorielles}. It was also found by Grinstead \cite{Grinstead:note}
in the case $k=2$, and it is related to some results of Grabiner
\cite{Grabiner:brownian}.

The idea of the argument is that, if we 
name the dimensions appearing in equation~\eqref{e:decomp},
$$f_\lambda = \dim R_\lambda \qquad 
d_\lambda = \dim V_\lambda \qquad
n_\lambda = f_\lambda d_\lambda,$$
the quantity $f_\lambda$ can be considered in the context of the
\emph{$k$-ballot problem}.  Suppose $N$ voters vote sequentially for an ordered
list of $k$ candidates.  In how many ways can they cast their votes so that the
$a$th candidate is never ahead of the $b$th candidate for $a>b$, and at the end
the $a$th candidate has $\lambda_a$ votes for every $a$?  Such a sequence of
votes is a \emph{ballot sequence} of shape $\lambda$ and there are $f_\lambda$
of them \cite[Prop. 7.10.3]{Stanley:vol2}. In this context of the RSK algorithm,
$f_\lambda$ is the number of standard tableaux of shape $\lambda$.  If the
$t$th entry of such a tableau is in row $a$, we can say that the $t$th voter
votes for candidate $a$.  This establishes a bijection between standard
tableaux and ballot sequences.  That $f_\lambda$ is the number of standard
tableaux of shape $\lambda$ can also be seen  directly from the representation
theory of $\gl(k)$:  The generalized Clebsch-Gordan rule states that
$$V \tensor V_\lambda = \bigoplus_{\lambda' = \lambda + \Box} V_{\lambda'},$$
where the sum is over shapes $\lambda'$ that are obtained from $\lambda$ by
adding a single box.  Thus the multiplicity of $V_{\lambda}$ in $V^{\tensor N}$
is the number of increasing chains of partitions from the empty partition to
the partition $\lambda$.  Such a chain is equivalent to a standard tableau of
$\lambda$ by assigning $t$ to a box if it appears at step $t$. 

In this formulation, the number $f_\lambda$ may be computed by the
\emph{reflection principle} \cite{GZ:walk}, which is also a disguised version
of the Weyl character formula \cite{GM:random}.  Recognizing the set of
partitions as a subset of $\Z^k$, there is an action of the symmetric group
$S_k$ on $\Z^k$ given by permuting coordinates. A chain of partitions is then a
lattice path in $\Z^k$ that happens to stay in the cone of partitions.  Here a
valid lattice path is one which increases one coordinate by one at each step.
Consider the partition
$$\delta = (k-1,k-2,k-3\ldots,2,1,0)$$
and let $m_\lambda$ for any $\lambda \in \Z^k$ be the number
of lattice paths from the origin to $\lambda$.  The reflection
principle shows that
\eq{e:reflect}{f_\lambda = \sum_{\sigma \in S_k}
    (-1)^\sigma m_{\lambda + \delta - \sigma(\delta)}}
Equation~\eqref{e:reflect} says that the number of ballot sequences from 0 to
$\lambda$ is the alternating sum of unrestricted lattice paths from a set of
image points of the form $\sigma(\delta) - \delta$ to $\lambda$. 
Figure~\ref{f:reflect} shows an example of the principle when $k=3$; in the
figure, partitions $\lambda$ are replaced by $\hlambda$ to obtain walks in a
2-dimensional lattice.

\begin{fullfigure}{f:reflect}
    {Paths from $\hlambda$ to $0$ and an image point in the weight lattice of
    $\su(3)$.}
\pspicture(-3.5,-2.165)(3.5,5.629)
\multirput(-3,5.196)(1,0){7}{\qdisk(0,0){1.5pt}}
\multirput(-2.5,4.33)(1,0){6}{\qdisk(0,0){1.5pt}}
\multirput(-2,3.464)(1,0){5}{\qdisk(0,0){1.5pt}}
\multirput(-1.5,2.598)(1,0){4}{\qdisk(0,0){1.5pt}}
\multirput(-1,1.732)(1,0){4}{\qdisk(0,0){1.5pt}}
\multirput(-1.5,.866)(1,0){5}{\qdisk(0,0){1.5pt}}
\multirput(-2,0)(1,0){5}{\qdisk(0,0){1.5pt}}
\multirput(-1.5,-.866)(1,0){4}{\qdisk(0,0){1.5pt}}
\multirput(-1,-1.732)(1,0){3}{\qdisk(0,0){1.5pt}}
\psline(-3.25,5.629)(1.25,-2.165)\psline(-2.5,0)(2.5,0)
\psline(-1.25,-2.165)(3.25,5.629)
\rput[b](.5,4.58){$\hlambda$}
\psline(.5,4.33)(-.5,4.33)(0,3.464)(.5,2.598)(-.5,2.598)(0,1.732)
\psline[linestyle=dashed](.5,4.33)(1,3.464)(1.5,2.598)
    (2,1.732)(2.5,.866)(1.5,.866)
\pscircle*[linecolor=white](0,1.732){1ex}
\pscircle*[linecolor=white](1.5,.866){1ex}
\rput(0,1.732){$+$}     \rput(1.5,.866){$-$}   \cput*(-1.5,.866){$-$}
\cput*(-1.5,-.866){$+$} \cput*(1.5,-.866){$+$} \cput*(0,-1.732){$-$}
\endpspicture
\end{fullfigure}

Since the numbers $m_\lambda$ are defined by walks with the same $k$ possible
steps at each time $t$, they can be approximated by the local central limit
theorem:
\eq{e:mlam}{m_\lambda \sim C e^{-k\hlambda^2/2N}}
Here and below we assume that $C$ is a constant depending only on $N$
and $k$, and we use the notation
$$\lambda^2 = \sum_{a=1}^k \lambda_a^2.$$
If the approximation~\eqref{e:mlam} were robust with respect to local finite
differences, then by the reflection principle it would give us an estimate for
$f_\lambda$.  If it were also robust under amplification by a polynomial
in $\lambda$, it would give us an estimate for
$$n_\lambda = f_\lambda d_\lambda,$$
since the Weyl dimension formula \cite[\S24.3]{Humphreys:gtm} says
that 
\eq{e:dim}{d_\lambda = \prod_{a>b} \frac{\lambda_a-\lambda_b + a-b}{a-b}}
is a polynomial in $\lambda$. Both of these refinements of the local central
limit theorem are true for arbitrary bounded lattice walks.
To state the theorem, we introduce a few definitions.

A \emph{finite difference operator} $D$ is a linear operator on functions
$$p:\R^k \to \R$$
defined by a finite sum
$$Dp(v) = \sum_t c_t f(v + v_t)$$
for some constants $\{c_t\}$ and some vectors $\{v_t\}$. The \emph{degree} $a$
of $D$ is the minimum degree of a polynomial $p$ such that $Dp \ne 0$. If $L
\subset \R^k$ is a lattice, the determinant
$$\det L = \Vol\ \R^k/L$$
is defined as the volume of the quotient space, or equivalently as the
determinant of a positive basis for $L$.

\begin{theorem} Let $X$ be a bounded, mean 0 random variable taking values in
$z+L$ for some lattice $L \subset \R^k$ and some vector $z \in \R^k$.  Assume
that $L$ is the thinnest such lattice for the given $X$.   Let
$$X' = \sum_{t=1}^N X^{(t)}$$
denote the sum of $N$ independent copies of $X$. Let
\begin{align*}
p(v) &= P[X' = v] \\
q(v) &= \frac{\det L}{(2\pi)^{k/2}\sqrt{\det \kappa}} e^{-\kappa^{-1}(v,v)/2N},
\end{align*}
where $\kappa$ is the covariance form of $X$.
Then for every finite difference operator $D$ of degree $a$ and for every
integer $b \ge 0$,
$$\lim_{N \to \infty} N^{(k-b+a)/2} |v|^b D(p-q)(v) = 0$$
uniformly for $v \in Nz+L$.
\label{th:lcentral}
\end{theorem}

Lawler~\cite[Th. 1.2.1]{Lawler:walks} proves a special case of 
Theorem~\ref{th:lcentral} in which $a$ and $b$ are 0 or 2, $X$ has the uniform
distribution among nearest-neighbor steps in $\Z^k$, and $D$ has a restricted
form.  However, the proof actually establishes Theorem~\ref{th:lcentral} in its
almost its full generality, requiring only that $b$ be even.   The conclusion
for $b$ odd follows by taking the geometric mean of the formulas for $b-1$ and
$b+1$.

Since $f_\lambda$ is given by the hook-length formula \cite[Cor.
7.21.6]{Stanley:vol2}, we can also prove Theorem~\ref{th:lcentral} in this
special case using Stirling's approximation \cite[\S4]{Johansson:plancherel}.
Such a special argument is analogous to the special argument for the
Laplace-de Moivre theorem, which is the simplest case of the usual central
limit theorem.  But it is not enough for our later generalization,
Theorem~\ref{th:general}.

\begin{corollary} If $\lambda \vdash N$, then 
$$\lim_{N \to \infty} N^{k/2} (k^{-N} n_\lambda - 
C\prod_{a<b} (\lambda_a-\lambda_b)^2 e^{-k \hlambda^2/2N}) = 0$$
uniformly in $\lambda$.
\label{co:main}
\end{corollary}
\begin{proof} We apply Theorem~\ref{th:lcentral}.  First, we change 
$\lambda$ from a subscript to an argument in certain quantities
that depend on it (and implicitly on $N$):
$$f(\hlambda) = f_\lambda\qquad m(\hlambda) = m_\lambda \qquad
n(\hlambda) = n_\lambda,$$
where dependence on $N$ is implicit in the notation.
Let $L = \Lambda$
be the set of all centered partitions $\hlambda$ (the weight lattice).
Define a finite difference operator $D$ by
$$Dp(\hlambda) = \sum_{\sigma \in S_k}
(-1)^\sigma p(\hlambda + \hdelta - \sigma(\hdelta))$$
so that
$$Dm(\hlambda) = f(\hlambda)$$
by equation~\eqref{e:reflect}.  The two important properties of the operator
$D$ are first, that it is antisymmetric under the Weyl group $S_k$ after
translation by $\hdelta$, and second, that it has degree $k(k-1)/2$.
The degree of $D$ follows from a factorization that appears
in proofs of the Weyl dimension formula \cite[\S24.3]{Humphreys:gtm}:
$$D = \prod_{\alpha \in \Phi_+} D_\alpha,$$
where $\Phi_+$ is the set of positive roots of $\su(k)$ and
$$D_\alpha(\hlambda) = p(\hlambda) - p(\hlambda+\alpha).$$
Each $D_\alpha$ has degree 1 and there are $k(k-1)/2$ of them.
Note that the only antisymmetric polynomial of degree $k(k-1)/2$
is
$$\Delta = \prod_{a<b} (\hlambda_a-\hlambda_b)
= \prod_{a<b} (\lambda_a-\lambda_b),$$

When $N$ is large,
$$De^{-k \hlambda^2/2N} \sim C \Delta e^{-k\hlambda^2/2N}$$
because in the limit $D$ becomes an antisymmetric differential operator of
degree $k(k-1)/2$.  When applied to a symmetric Gaussian, it produces an
anti-symmetric polynomial factor of degree $k(k-1)/2$. The polynomial $\Delta$
is the only choice for this factor up to scale. Thus the operator $D$ explains
one factor of $\Delta$ in the statement of the corollary.  The other factor is
given by $d_\lambda$, which is also proportional to $\Delta$ in the limit as
$\lambda \vdash N$ goes to $\infty$. Theorem~\ref{th:lcentral} then establishes
the stated approximation for $n_\lambda$, where $D$ and $d_\lambda$ each
contribute a factor of $\Delta$.
\end{proof}

Corollary~\ref{co:main} is evidently the precise statement of
Theorem~\ref{th:main}.

\section{Generalizations}
\label{s:general}

The first way that we can generalize Theorem~\ref{th:main} that is that we can
replace the representation $V$ of $\su(k)$ by some other finite-dimensional
representation $W$.  In general a tensor power of such a representation
decomposes as
\eq{e:wdecomp}{W^{\tensor N}
    \cong \bigoplus_{\lambda\vdash 0} T_{\lambda,N} \tensor V_{\lambda},}
where each $T_{\lambda,N}$ is a vector space on which $\su(k)$ acts trivially. 
(In this generality it does not make sense to make $\lambda$ a partition of $N$
or any other particular integer, so we take it to be a highest weight, or a
partition of 0.)  The space $T_{\lambda,N}$ is a representation of the
symmetric group $S_N$, but it is not usually irreducible, not even when $W$
is.  Assuming the a state on $W^{\tensor N}$ which is invariant under the
action of $\su(k)$, we may as before use equation~\eqref{e:wdecomp} to define a
quantum random variable $\lambda_{\QM}$.

It is less trivial to define a classical counterpart $\lambda_{\RSK}$, or even
the space of words on which it is defined.  If $W$ is irreducible, we
can model it as a summand of $V^{\tensor \ell}$ for
some $\ell$.  More precisely, we choose a partition $\mu \vdash \ell$
such that $W \cong V_{\hmu}$, and we choose a specific standard 
tableau $Q_W$ of shape $\mu$.  Then the set $S \subset [k]^\ell$ of words 
with recording tableau $Q_W$ indexes a basis of $W$.  The set
$S$ can be interpreted as a ``syllabic alphabet'', in the sense
that a word $w$ of length $N$ over the alphabet $S$ is simultaneously a
word of length $\ell N$ over the alphabet $[k]$.  Remarkably, the RSK algorithm
is compatible with this dual interpretation:  If we define
the shape $\lambda$ of $w \in S^N$ by spelling it out in $[k]^{\ell N}$
and taking the usual shape,  then once again
$$\dim T_{\lambda,N} \tensor V_{\lambda} = n_\lambda,$$
where $n_\lambda$ is the number of words $w$ with shape $\lambda$. (One way to
argue this fact is with the theory of Littelmann paths; see below.  Syllabic
expansion corresponds to concatenation of paths.)

For example, if $k = 2$, then $V^{\tensor 2}$ has a summand $W = V_2$
isomorphic to the adjoint representation of $\su(2)$. As it happens, this
summand occurs only once.  If we take left and right parentheses
$\{\:)\:,\:(\:\}$ as the alphabet for the basis of $V$ rather than $\{1,2\}$,
then the first component $\hlambda_1$ of the centered shape of a string is half
the number of unmatched parentheses.  For example, $\hlambda_1 = 1$ for the
string
$$\pspicture(-.5,0)(3.5,1.75)\begin{large}
\rput[b](0,0){\rnode{u1}{$)$}}\rput[b]( .5,0){\rnode{a1}{$($}}
\rput[b](1,0){\rnode{b1}{$)$}}\rput[b](1.5,0){\rnode{a2}{$($}}
\rput[b](2,0){\rnode{a3}{$($}}\rput[b](2.5,0){\rnode{b3}{$)$}}
\rput[b](3,0){\rnode{b2}{$)$}}
\nccurve[nodesep=3pt,angle=90,ncurv=1.5]{a1}{b1}
\nccurve[nodesep=3pt,angle=90,ncurv=1.5]{a2}{b2}
\nccurve[nodesep=3pt,angle=90,ncurv=1.5]{a3}{b3}
\end{large}\endpspicture$$
The alphabet $S$ for the representation $W$ is the set of three pairs of
parentheses $\{\:((\:,\:)(\:,\:))\:\}$ other than the two that match each
other.   If we rename this alphabet $\{\:\langle\:,\:|\:,\:\rangle\:\}$,
then one can check that the only words that match
completely are those that form nested complete ``bra-kets'':
$$\pspicture(-.5,0)(4,.75)\begin{large}
\rput[b](0,0){$\langle$}\rput[b](.5,0){$|$}
\rput[b](1,0){$|$}\rput[b](1.5,0){$\langle$}
\rput[b](2,0){$|$}\rput[b](2.5,0){$\rangle$}
\rput[b](3,0){$|$}\rput[b](3.5,0){$\rangle$}
\psline(-.15,.1)(-.15,-.3)(3.65,-.3)(3.65,.1)
\psline(1.35,.1)(1.35,-.15)(2.65,-.15)(2.65,.1)
\end{large}\endpspicture$$
A general string will have a maximal substring of this form, as well as
fragments consisting of unmatched ``bras'', unmatched ``kets'', and
unbracketed separators:
$$\mbox{\large $\langle\;|\;|\;| \qquad |\;|\;|\;\rangle \qquad |\;|\;|\;|$}$$
The statistic $\hlambda_1$ is then the number of these fragments.
For example, $\hlambda_1 = 2$ for the string
$$\pspicture(-.5,0)(4,.75)\begin{large}
\rput[b](0,0){$|$}\rput[b](.5,0){$|$}\rput[b](1,0){$\langle$}
\rput[b](1.5,0){$|$}\rput[b](2,0){$|$}\rput[b](2.5,0){$\rangle$}
\rput[b](3,0){$\langle$}\rput[b](3.5,0){$|$}\rput[b](4,0){$|$}
\psline(-.15,-.15)(.65,-.15)
\psline(.85,.1)(.85,-.15)(2.65,-.15)(2.65,.1)
\psline(2.85,.1)(2.85,-.15)(4.15,-.15)
\end{large}\endpspicture$$
since there are four unmatched parentheses if the bra-kets 
are expanded into parentheses:
$$\pspicture(-.5,0)(9,1.75)\begin{large}
\rput[b](0,0){\rnode{u1}{$)$}}\rput[b](0.4,0){\rnode{a1}{$($}}
\rput[b](1,0){\rnode{b1}{$)$}}\rput[b](1.4,0){\rnode{u2}{$($}}
\rput[b](2,0){\rnode{a2}{$($}}\rput[b](2.4,0){\rnode{a3}{$($}}
\rput[b](3,0){\rnode{b3}{$)$}}\rput[b](3.4,0){\rnode{a4}{$($}}
\rput[b](4,0){\rnode{b4}{$)$}}\rput[b](4.4,0){\rnode{a5}{$($}}
\rput[b](5,0){\rnode{b5}{$)$}}\rput[b](5.4,0){\rnode{b2}{$)$}}
\rput[b](6,0){\rnode{u3}{$($}}\rput[b](6.4,0){\rnode{a6}{$($}}
\rput[b](7,0){\rnode{b6}{$)$}}\rput[b](7.4,0){\rnode{a7}{$($}}
\rput[b](8,0){\rnode{b7}{$)$}}\rput[b](8.4,0){\rnode{u4}{$($}}
\psline(-.15,-.15)(1.6,-.15)
\psline(1.8,.1)(1.8,-.15)(5.6,-.15)(5.6,.1)
\psline(5.8,.1)(5.8,-.15)(8.6,-.15)
\nccurve[nodesep=3pt,angle=90,ncurv=1.5]{a1}{b1}
\nccurve[nodesep=3pt,angle=90,ncurv=.75]{a2}{b2}
\nccurve[nodesep=3pt,angle=90,ncurv=1.5]{a3}{b3}
\nccurve[nodesep=3pt,angle=90,ncurv=1.5]{a4}{b4}
\nccurve[nodesep=3pt,angle=90,ncurv=1.5]{a5}{b5}
\nccurve[nodesep=3pt,angle=90,ncurv=1.5]{a6}{b6}
\nccurve[nodesep=3pt,angle=90,ncurv=1.5]{a7}{b7}
\end{large}\endpspicture$$
More generally still, we can let $W$ be any non-trivial, finite-dimensional,
unitary representation of any compact simple Lie algebra $\mg$.  (We say that a
Lie algebra is \emph{compact} if it integrates to a compact Lie group.)  Once
again there is a direct sum decomposition
$$W^{\tensor N}
   \cong \bigoplus_{\lambda \in \Lambda} T_{\lambda,N} \tensor V_{\lambda},$$
where $\Lambda$ is the weight lattice of $\mg$ and $V_{\lambda}$
is the irreducible representation of highest weight $\lambda$.  As before
this decomposition defines a quantum random variable $\lambda_\QM$
if we assume the tracial state $\rho$ on $W$.

If $W$ is irreducible, the theory of Littelmann paths then provides a
satisfactory combinatorial counterpart $\lambda_\LP$ with the same distribution
as $\lambda_\QM$ \cite{Littelmann:paths}.  If $W \cong V_\mu$, then we can apply
the Littelmann lowering operators to some fixed dominant path $p_\mu$ from the
origin to $\mu$. There is a natural bijection between the resulting set of
paths $P(W)$ and a basis of $V_\mu$.  Moreover, a word $w$ in $P_\mu$ forms a
longer path $\gamma(w)$ given by concatenating letters.  If we apply Littelmann
raising operators as many times as possible to $\gamma(w)$, the result is a
highest weight $\lambda$.  Assuming the uniform distribution on $P_\mu$, this
weight defines a random variable random variable $\lambda_\LP = \lambda$. Note
that $\lambda_\LP$ depends on $p_\mu$, although its distribution does not.

Although the abstract setting of Littelmann paths looks quite different from
the Robinson-Schensted-Knuth algorithm, it is in fact a strict generalization
\cite{Leeuwen:private}. Briefly, if $W = V$ is the defining representation of
$\su(k)$ and $p_\mu$ is a straight line segment, then every element of $P_\mu$
is a straight line segment, and these segments are naturally enumerated by the
integers $1,\ldots,k$.  The highest weight $\lambda_\LP$ of a word $w$
coincides with the centered shape of the tableau produced by a \emph{dual} RSK
algorithm defined using column insertion.  (The standard RSK algorithm
uses row insertion.)  By one of the symmetries of the RSK algorithm \cite[Cor.
A1.2.11]{Stanley:vol2}, this shape is the same as the row-insertion shape of
the reverse word $w^*$.  Thus
$$\lambda_\LP(w) = \hlambda_\RSK(w^*) \qquad \lambda_\LP \doteq \hlambda_\RSK.$$

Finally, if $W$ is not irreducible, then we can choose a separate alphabet for
each summand in a direct-sum decomposition. For example, if
$$W \cong V_\mu \oplus V_\mu$$
for some dominant weight $\mu$, then we can let the alphabet be the disjoint
union of two copies of the same alphabet $P_\mu$, a ``red'' copy and a ``blue''
copy.  The fundamental properties of Littelmann paths imply that in all cases,
$\lambda_\QM \doteq \lambda_\LP$, where
$$P[\lambda_\QM = \lambda] = \frac{\dim T_{\lambda,N} \tensor V_{\lambda}}
{(\dim W)^N}.$$

What is the counterpart to $\lambda_\GUE$ in this context?
Taking Section~\ref{s:quantum} as a guide, each element $A \in i \cdot \mg$
defines a real-valued random variable on the quantum space $(W,\rho)$.
These variables have the covariance form
$$\kappa_W(A,B) = \rho(AB) = \frac{\Tr_W(AB)}{\dim W}.$$
All together they can be taken as 1-dimensional projections  of a quantum
random variable $x \in i\cdot \mg^*$. Instead of the spectrum, we can consider
the orbit of $ix$ under the co-adjoint action of $\mg$ on $i\cdot \mg^*$. By
standard representation theory, $ix$ is conjugate to a unique weight
$$\lambda \in \mC \subset \mh^* \subset \mg^*,$$
where $\mC$ is a Weyl chamber in $\mh^*$, the dual space to a Cartan subalgebra
of $\mg$. If we assume a Gaussian distribution $\mu_W$ on $i\cdot\mg^*$ with
covariance matrix $\kappa_W$, the corresponding distribution for the weight
$\lambda$ is
$$e_W(\lambda) d\lambda = C\prod_{\alpha \in \Phi_+} \alpha(\lambda)^2
e^{-\kappa_W^{-1}(\lambda,\lambda)/2} d\lambda,$$
where as before $\Phi_+$ is the set of positive roots in $\mg$. This
distribution can be derived in the same way as equation~\eqref{e:mehta}. If
$\mg = \sp(2n)$ (the compact form of $\sp(2n,\C)$), it is the Gaussian
symplectic ensemble (GSE). But if $\mg = \so(n)$, it is the Gaussian
antisymmetric ensemble (GAE).  Since symmetric matrices do not form a Lie
algebra, the Gaussian orthogonal ensemble (GOE) would require some yet more
general model.

Finally we can state the general theorem.

\begin{theorem} Let $W$ be a non-trivial, finite-dimensional, unitary
representation of a compact simple Lie algebra $\mg$ of rank $r$.  Let
$n_\lambda$ be the dimension of the isotypic summand of $W^{\tensor N}$ of
highest weight $\lambda$.  Then
$$\lim_{N \to \infty}
r^{N/2} \bigl(\frac{n_\lambda}{(\dim W)^N} - Ce_W(\frac{\lambda}{\sqrt{N}})\bigr) = 0$$
uniformly in $\lambda$.
\label{th:general}
\end{theorem}

The arguments of Sections~\ref{s:quantum} and \ref{s:local} both generalize in
a straightforward way to proofs of Theorem~\ref{th:general}.  As before,
Section~\ref{s:quantum} establishes a weak version of it at a rigorous level.
We also comment that the tautological matrix $M$ of Section~\ref{s:quantum}
should be replaced by a certain  $i\cdot\mg^*$-valued measurement operator
$$M \in \mg^* \tensor \mg$$
acting on $W^{\tensor N}$.  As a tensor in $\mg^* \tensor \mg$, $M$ is again
tautological; it comes from the identity linear transformation from $\mg$ to
itself.

\begin{remark} The Lie algebra picture of Theorem~\ref{th:general} suggests
another interpretation which is dual to that of  Section~\ref{s:quantum}, and
in another sense dual to that of Section~\ref{s:local}.  If $G$ is a compact,
simple Lie group with Lie algebra $\mg$ and $W$ is a unitary representation of
$G$, then the absolute value of the  character $\chi_W$ of $W$ has a local
maximum at $1 \in G$. When $N$ is large, the character $\chi_W^N$ of
$W^{\tensor N}$ is approximately a Gaussian in a neighborhood of 1. If we
inflate $G$ by a factor of $\sqrt{N}$, it converges to $\mg$, and multiplication
on $G$ converges to addition on $\mg$. The character $\chi_W^N$ converges to a
limit on $\mg$, namely the Fourier transform of the Gaussian distribution
$\mu_W$ on $\mg^*$ defined above.  This intermediate picture led the author
from Section~\ref{s:local} to Section~\ref{s:quantum}.
\end{remark}

\subsection{Things out of reach}

When $\mg = \su(k)$, Theorem~\ref{th:general} can be interpreted as a limit
distribution result for the shape $\lambda_\RSK$ of words with various
interesting distributions.  For example, if each letter of the alphabet for the
representation $V_2$ of $\su(2)$ is expanded into a pair of letters in the
alphabet $\{\up,\down\}$, then the distribution $\rho$ on expanded words is
determined by its correlations for \emph{digraphs} (adjacent pairs of
letters).  In this case the correlation between the $t$th and $t+1$st letter
depends on whether $t$ is odd or even.  But random words associated with
representations such as $V_2 \oplus V_3$ do not exhibit such irregularities.

Especially when $k=2$, these distributions resemble distributions given by
\emph{doubly stochastic Markov chains}.  In other words, the first letter $w_1$
of a random word $w \in [k]^N$ has the uniform distribution. Each subsequent
letter depends on the immediate predecessor (but not on earlier letters)
according to a Markov matrix $M$:
$$P[w_{t+1} = a \st w_t = b] = M_{ab}.$$
Here $M$ is chosen so that every letter has the uniform distribution
if the first one does.

What is the asymptotic distribution of the shape $\lambda_\RSK$
of a random word $w \in [k]^N$ generated by a doubly stochastic Markov
matrix $M$?  Non-rigorously we expect it to have the form
$$C P(\hlambda) e^{-k\hlambda^2 / 2vN}.$$
Here $P$ is some polynomial (or at least some function which is asymptotically
polynomial) and $v$ is the \emph{variance per letter} of $w$.  The variance $v$
is defined by the formula
$$v = \sum_{t = -\infty}^{\infty} \frac{kP[w_0 = w_t]-1}{k-1}$$
using a bi-infinite word $w$ generated by $M$.

\begin{figure}[htb]\begin{center}
\psfrag{l1}{\small $A$}
\psfrag{l2}{\small $F$}
\psfrag{l3}{\small $C_+$}
\psfrag{l4}{\small $C_-$}
\psfrag{lx}{\small $\hlambda_1$}
\includegraphics*[65pt,50pt][310pt,239pt]{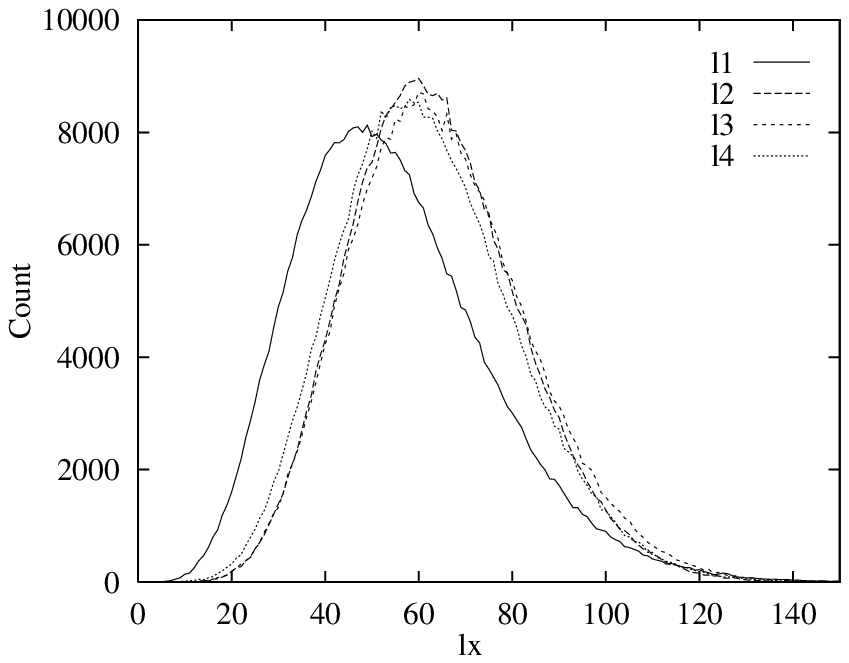}
\caption{Distribution of $\hlambda_1$, the centered length of the longest
increasing subsequence, for words generated by four different Markov chains,
400,000 trials each.}
\label{f:plot}
\end{center}\end{figure}

We have conducted computer experiments with different choices of $M$
with 2- and 3-letter alphabets \cite{Kuperberg:words.c}.  Figure~\ref{f:plot}
shows the distribution of $\hlambda_1$ for 400,000 words generated
by each of the following four Markov matrices $M$:
\begin{align*}
A   &= \frac14\matrix{ccc}{1&2&1 \\ 3&1&0 \\ 0&1&3} &
F   &= \frac13\matrix{ccc}{1&1&1 \\ 1&1&1 \\ 1&1&1} \\[1ex]
C_+ &= \frac14\matrix{ccc}{3&0&1 \\ 1&3&0 \\ 0&1&3} &
C_- &= \frac14\matrix{ccc}{3&1&0 \\ 0&3&1 \\ 1&0&3}
\end{align*}
The lengths of the words are 1620, 3420, 1140, and 1140 in the four respective
cases. These lengths were chosen so that the four types of words would have the
same total variance (ignoring boundary effects).  The experiments indicate that
the distribution of $\hlambda_1$ (the centered length of the longest weakly
increasing subsequence) in the asymmetric distribution $A$ is genuinely
different from the referent uniform distribution $F$.  The lower median value
of $\hlambda_1$ in this case does not disappear as the words grow longer.  It
also cannot be explained as a maladjusted variance, because at the other end
the tail of $A$ eventually overtakes the tail of $F$.  On the other hand, the
distribution for the cyclic Markov chains $C_+$ and $C_-$ do appear to converge
to the distribution for $F$.  Their symmetry implies that the longest
weakly increasing subsequence sees the same fluctuations in the transition $1
\to 2$ as it does for the transition $2 \to 3$, which is apparently
enough to produce the same distribution.

\begin{conjecture} Let $M$ be an indecomposable, doubly stochastic
matrix such that
$$M_{a,b} = M_{a+1,b+1},$$
where $k+1 \equiv 1$.  Then the distribution of the shape
of a word of length $N$ generated by $M$ converges locally
to the distribution of the spectrum of a traceless $k \times k$
GUE matrix.
\label{c:cyclic}
\end{conjecture}

Conjecture~\ref{c:cyclic} can be generalized further by considering
other cyclically symmetric, translation-invariant measures on words
whose correlations decay sufficiently quickly.

\begin{problem} Let $V = \C^k$ and let $\rho$ be the state on
$V^{\tensor n}$ extending a distribution on $[k]^N$ generated
by a doubly stochastic Markov chain.  What is the limiting distribution
of $\lambda_\QM$?
\label{p:qm}
\end{problem}

Problem~\ref{p:qm} is really a statistical mechanics question concerning a
quantum spin chain with certain nearest-neighbor interactions. It cannot be
stated in terms of $\lambda_\RSK$ because there is no reason to expect that
$\hlambda_\RSK \doteq \lambda_\QM$ in this generality. Yet more generally, we
can ask about the behavior of $\lambda_\QM$ for an arbitrary nearest-neighbor
interaction that produces the tracial state when restricted to a single site. 

\begin{problem} What is the distribution of the longest
weakly increasing circular subword of a circular word $w \in [k]^N$?
\label{p:circular}
\end{problem}

In Problem~\ref{p:circular}, we assume that both $[k]$ and $[N]$  are
circularly ordered.  We do not know if there is a suitable circular
generalization of the RSK algorithm.

\subsection{Infinite matrices}

The most interesting case of Theorem~\ref{th:main} to consider (indeed the case
that motivated the result) is the limit $k \to \infty$. We have no firm 
results about this limit, but we can propose a model of it that may be
important.  Our model might be related to the semi-infinite
wedge space model of Okounkov \cite{Okounkov:wedge}.

\begin{fullfigure}{f:II1}{The hyperfinite $\mathrm{II}_1$ factor $\cM$
as a matrix algebra over itself}
\begin{Large}$\cM\;\cong$\;\end{Large}
\pspicture[.45](0,0)(6,6)
\multirput(0,0)(1.2,0){6}{\psline[linecolor=gray40](0,0)(0,6)}
\multirput(0,0)(0,1.2){6}{\psline[linecolor=gray40](0,0)(6,0)}
\multirput(.6,0)(1.2,0){5}{\lightgray\multirput(0,.6)(0,1.2){5}{$\cM$}}
\multirput(0,0)(2,0){4}{\psline(0,0)(0,6)}
\multirput(0,0)(0,2){4}{\psline(0,0)(6,0)}
\multirput(1,0)(2,0){3}{\multirput(0,1)(0,2){3}{\large $\cM$}}

\endpspicture
\end{fullfigure}

Consider the Hilbert space $\cH = L^2([0,1])$.  For every 
$k$, the matrix algebra $\cM_k$ acts on $\cH$ by taking
$$E_{ab}(f)(x) = f(x+\frac{b-a}k)$$
if $f$ is supported on $[\frac{b-1}k,\frac{b}k]$, and
$$E_{ab}(f)(x) = 0$$
if $f$ vanishes on $[\frac{b-1}k,\frac{b}k]$.
The weak-operator closure of all of these algebra actions is a von Neumann
algebra $\cM$ called the \emph{hyperfinite $\mathrm{II}_1$ factor}
\cite[\S12.2]{KR:vol2}.  For every $k$, $\cM$ is a $k \times k$
matrix algebra over itself (Figure~\ref{f:II1}).  Thus it 
could be generally important in random matrix theory.

In this case we are interested in the Lie algebra structure of $\cM$ (in
addition to its topologies), making $\cM$ an infinite-dimensional analogue of
$\gl(k,\C)$. We are also interested in the tracial state $\rho$, defined as the
continuous extension of the tracial state on each $\cM_k$.  It is a normal
state.  It is also a model of the uniform measure on the interval $[0,1]$. 
Finally we define $\hat{\cM}$ to be the kernel of $\rho$, analogous to the space
$\sl(k,\C)$ of traceless matrices.

As a Lie algebra, $\cM$ acts on $\cH^{\tensor N}$ for every $N$.  This action
commutes with the action of the symmetric group $S_N$ on $\cH^{\tensor N}$
given by permuting tensor factors.  There is a direct-sum decomposition
\eq{e:infsum}{\cH^{\tensor N} \cong
\bigoplus_{\lambda \vdash N} \cH_\lambda,}
where $\cH_\lambda$ is, as a representation of $S_N$, the isotypic summand of
type $R_\lambda$.  Each of these representations has a measure-theoretic
dimension defined using $\rho^{\tensor N}$ on $\cH^{\tensor N}$ (in which $\cH$
embeds by the usual Leibniz rule for Lie algebra actions):
$$\dim \cH^{\tensor N} = 1 \qquad \dim \cH_\lambda = \frac{f_\lambda^2}{N!}.$$
Thus equation~\eqref{e:infsum} is a quantum statistics model for the Plancherel
measure on the symmetric group. For each $N$, it defines a quantum random
variable $\lambda_\QM \vdash N$.

The state $\rho^{\tensor N}$ also expresses the uniform measure on $[0,1]^N$,
\ie, the process of choosing a ``word'' of $N$ random points in the unit
interval.   The usual RSK algorithm is defined for such words.  Since the
letter of the word are distinct almost surely, and since the RSK algorithm
depends only on the order of the letters and not their values, it defines a
random variable $\lambda_\RSK$ equivalent to the shape of a random permutation.
Its distribution is also the Plancherel measure.

By the quantum central limit theorem, the state $\rho^{\tensor N}$ should
produce a Gaussian measure on $\hat{\cM}$ in the limit $N \to \infty$.  So
should the GUE measure on $\sl(k,\C)$ in the limit $k \to \infty$.  The
relation between these two limits could shed light on the Vershik-Kerov limit
for Plancherel measure \cite{VK:limit} and the Wigner semicircle for the
spectrum of a random matrix \cite{Mehta:matrices}.  The quantum central limit
theorem might also predict the distribution of the deviation from a semicircle,
at least to first order in $N$.


\providecommand{\bysame}{\leavevmode\hbox to3em{\hrulefill}\thinspace}

\end{document}